\documentclass{article}
\usepackage{amssymb}
\begin{document}

\newcommand{\dsp}{\displaystyle}
\newcommand{\eps}{\varepsilon}
\newcommand{\tends}{\rightarrow}
\newcommand{\kraj}{\ $\Box$}

\newcommand{\R}{\bold{R}}
\newcommand{\cF}{{\mathcal F}}
\newcommand{\cS}{{\mathcal S}}
\newcommand{\cB}{{\mathcal B}}

\newcommand{\di}{{\: \rm d}} 
\newcommand{\E}{{\rm E\, }} 


\newcommand{\be}{\begin{equation}}
\newcommand{\ee}{\end{equation}}

\newcommand{\beqn}{\begin{eqnarray*} }
\newcommand{\eeqn}{\end{eqnarray*}}


\newcounter{secti}
\newcounter{item}[secti]
\renewcommand{\theitem}{\thesecti.\arabic{item}}

\newcommand{\sect}[1]{\medskip \begin{center} \refstepcounter{secti}
{\sc \thesecti.\  #1 }\end{center} \medskip}

\newcommand{\qitem}[2]{\refstepcounter{item}
  {\bf \theitem\  #1.}\  {#2 }}

\newcommand{\qitems}[2]{\refstepcounter{item}
  {\bf \theitem\  #1}\  {#2 }}     

\newcommand{\eitem}{\medskip}


\begin{center}
{\Large\bf  An axiomatic integral and a multivariate  mean value theorem
}

\vspace{3em}

{\sc Milan Merkle}

\bigskip

\parbox{25cc}{\small
{\bf Abstract. }{\small In order to investigate minimal sufficient conditions for an abstract integral to belong to the convex hull of
the integrand, we propose a system of axioms  under which it happens.
 If the integrand is a continuous $\R^n$-valued function
 over a path connected topological space, we prove that any such integral can be represented as a convex combination of
 values of the integrand in at most $n$ points, which yields an ultimate multivariate mean value theorem.

}

\medskip



}

\end{center}

\medskip

\sect{Introduction and motivation}

\medskip

The basic integral mean value theorem states that for a function $X$ which is continuous on the interval $[a,b]$,
there exists a point $t\in (a,b)$ such
that
\be
\label{bastc}
\frac{1}{b-a} \int_a^b X(s) \di s =  X(t)
\ee
To show that the assumption of continuity is crucial for validity of this theorem, we can take the interval $[-1,1]$ and define $X(s)=-1$ for $s\in [-1,0)$
 and
$X(s)=1$ for $s\in [0,1]$. Hence here we do not have a single point $t\in (-1,1)$ for (\ref{bastc}) to be satisfied. However, we can achieve a similar
result with a convex combination of values of $f$ in two points:

\be
\label{bastd}
\frac{1}{b-a} \int_a^b X(s) \di s = \frac{1}{2}X(t_1) + \frac{1}{2}X(t_2),\qquad t_1\in (-1,0), \ t_2\in (0,1) .
\ee

It turns out that the difference of one extra point for non-continuous functions remains in much more general case and for a very broad class of integrals
in higher dimensions. This is the topic of this article.

In multivariate case, with  $X\in \R^n$, $n\geq 1$, there is an old and seemingly forgotten
 mean value theorem by Kowalevski \cite{kowal1,kowal2} as follows.

{\bf Theorem A\cite{kowal1}. }{\em Let $x_1,\ldots , x_n$ be
continuous functions in a variable $t\in [a,b]$. There exist real
numbers $t_1,\ldots, t_n$ in $[a,b]$ and  non-negative numbers
$\lambda_1,\ldots,\lambda_n$, with $\sum_{i=1}^n \lambda_i=b-a$,
such that
\[ \int_a^b X_k(t)\di t = \lambda_1X_k(t_1)+\cdots +\lambda_n X_k(t_n),\quad \mbox{for each $k=1,2,\ldots n$ },\]
or in more compact notation,
\be
\label{kowal1}
\int_a^b X(s)\di s = \lambda_1X(t_1)+\cdots +\lambda_n X(t_n),
\ee
where $X=(X_1,\ldots,X_n)$.
}

 A recent generalization of Theorem A  is
proved in \cite{kowalbm} using the following  modification of classical Carath\'eodory's convex hull theorem.

{\bf Theorem B\cite{kowalbm}.} {\em Let $C: s\mapsto X(s)$, $s\in I$, be a continuous curve in $\R^n$, where $I\subset \R$ is an interval,
and let $K$ be the convex hull of the  curve $C$.
Then each $v\in K$ can be represented as a convex combination of $n$ or fewer points of the curve $C$.}

In this way, a more general mean value theorem for $n$-dimensional functions
 is obtained directly from the fact that the normalized  integral should belong to the convex hull of the set of values of the
integrand. This is proved in the following theorem for Lebesgue integral with a probability measure.

{\bf Theorem C\cite[Lemma 1]{kowalbm}.}{\em Let $(S,\cF,\mu)$ be a probability space, and let $X_i: S\rightarrow \R$ , $i=1,\ldots, n$, be $\mu$-integrable
functions. Let $X(t)=(X_1(t),\ldots,X_n(t))$ for every $t\in S$. Then $\int_S X(t)\di \mu(t)\in \R^n$ is in the convex hull of the set
$X(S)=\{ X(t)\; |\; t\in S\}\subset \R^n$.}

Finally, the main result of \cite{kowalbm} reads

{\bf Theorem D \cite[Theorem 2]{kowalbm}.} {\em For an interval $I\subseteq \R$, let $\mu$ be a finite positive measure on
the Borel sigma-field of $I$.
 Let $X_k$, $k=1,\ldots , n$, $n\geq 1$, be continuous functions on $I$, integrable on $I$
with respect to the measure $\mu$. Then there exist points $t_1,\ldots, t_n$ in $I$,
and non-negative numbers $\lambda_1,
\ldots, \lambda_n$, with $\sum_{i=1}^n \lambda_i = \mu(I)$, such that }

\be
\label{finres}
\int_{I} X_k(s)\di \mu (s) = \sum_{i=1}^n \lambda_i X_k(t_i),\qquad k=1,\ldots, n.
\ee

Let us note that without the continuity assumption we still may use the Carath\'eodory's convex hull theorem which would yield
(\ref{finres}) with $n+1$ points $t_i$ and the same number of $\lambda_i$'s, which shows that the example at the beginning of the
text well describes the general situation in $\R^n$.

In this paper we give a more general theorem of the type (\ref{finres}), tracing the steps  of the proof in \cite{kowalbm}
 in a much more general context.
In Section 3 we show that a result like (\ref{finres}) holds if the integral over $I$ is replaced with a general linear functional on some function
space, under a system of  axioms, whereas the interval $I$ can be replaced by a topological space which is path connected.

To reach this goal, we need to extend the Theorem
C.  In section 2 we show that Theorem C holds for any linear functional which satisfies a condition slightly stronger than positivity, and
where functions $X_k$ are defined over an arbitrary nonempty set.

Applications of such a very general mean value theorem are numerous, and we are not discussing particular applications in this paper. Let us just mention
that as shown in \cite{kowalbm}, the theorems of this type can be considered as an aid to construct quadrature rules or their approximative versions,
see also a recent paper \cite{numfred13} for another application related to integrals. Another advantage of the approach presented in this paper is that
the results are widely applicable to different kinds of integrals treated as linear functionals over some space of functions.

\sect{Axioms and their consequences}

We start with an arbitrary nonempty set $S$ with an algebra $\cF$  (may be a sigma algebra as well) of its subsets. Therefore,
$\cF$ contains $S$ and if a set $A$ is result of finitely many set operations over sets in $\cF$, then $A\in \cF$.

Let  $\cS $ be a family  of functions $X:S\rightarrow \R$ which satisfies the following conditions:

\begin{itemize}
\item[C1] If $X_1,X_2 \in \cS$ then $aX_1+bX_2\in \cS$ for all $a,b \in \R$,
\item[C2] For  $B\in \cF$, the indicator function $I_B(\cdot)$ belongs to $\cS$,
\item[C3] For $X\in \cS$ and any interval $J$,  the set $\{s\in S:\ X(s)\in J\}$ is in  $\cF$.
\item[C4] For $X\in \cS$ such that $X(s)\geq 0$ for all $s\in S$, it holds that  $X\cdot I_{X\in J}\in \cS$  for any interval $J$.
\end{itemize}

Note that from C1 and C2 it follows that all constants are in $\cS$. Let us also note  that functions in $\cS$ are not assumed to be bounded.

Let $\E$ be a functional defined on $\cS$ and taking values in $\R$ such that
the following axioms hold.

\medskip

 \begin{enumerate}
 \item[{\bf A1}]  $\E c = c$ for any constant $c$;
 \item[{\bf A2}] $\dsp\E \left( \sum_{i=1}^m \alpha_i X_i \right) = \sum_{i=1}^m \alpha_i \E (X_i)$ for $a_j\in \R$ and  for  $X_i \in \cS$, $i=1,\ldots, m$
 \end{enumerate}

Now we may define a set function $P$ as

\be
\label{aprob0}
P(B) = \E (I_{B}(\cdot)),\qquad B\in \cF,
\ee

and consider yet another condition related to $P$:

\begin{itemize}
\item[C5] If $X\in \cS$ and $P(N)=0$, then $X\cdot I_N\in \cS$, and $\E (X\cdot I_N)=0$.
\end{itemize}

The last axiom that we propose is

\begin{enumerate}
  \item[{\bf A3}]  For $X\in \cS$, if  $\E X =0$ then either  $P(X=0)=1$  or there exist $s_1,s_2\in S$ such that $X(s_1)X(s_2)<0$,
   \end{enumerate}

or equivalently (see Lemma \ref{equ3} below)

\begin{enumerate}
  \item[{\bf A3$'$}]  For $X\in \cS$, if $X(s)\geq 0$ for all $s\in S$ and $\E X =0$, then $P(X=0)=1$.
   \end{enumerate}

Finally we  extend the functional $\E$ to act on functions with values in $\R^n$. Let  $X=(X_1,\ldots X_n)$ be a function from $S$ to $\R^n$,
where we assume that $X_i$, $i=1,\ldots,n$ satisfy axioms and conditions above, then we define
\[ \E (X) := (\E X_1, \ldots , \E X_n). \]

The central result of this section is Theorem \ref{ich}, where we show that under A1-A2-A3, $\E (X)$ belongs to the convex hull of $X(S)$.

 A similar  axiomatic  approach is applied in  Daniell's integral, and there are other axiomatic systems in the literature for different purposes like
 in \cite{gill93} for Riemann integrals in connection to evaluation the length of a curve,  general means in \cite{ostrowceb},
finitely additive probabilities (FAPs) in \cite{ssk2008} and applications in a recent article \cite{ssk2014}.
The system of axioms applied in this article differs from others in the literature in the conditions that allow non-absolute integrals, as well as in
 axiom 3, which is  slightly  stronger condition then
usual positivity. The reason of introducing this system of axioms is that
it  provides conditions under which   $\E X$ belongs to the convex hull of $X(S)$  (theorem \ref{ich}), independently of the kind
of integrals that is considered.

Now we are going to derive some additional properties as consequences from the axioms.

\qitem{Lemma}{\label{aprob} \em Under system of axioms A1-A2-A3 or A1-A2-A3$'$, assuming also conditions $C_1$-$C_2$,
the set function $P$ defined on $\cF$  with (\ref{aprob0})  is  a finitely additive probability on $(S,\cF)$.}

{\bf Proof. } Since $I_S (s) =1$ for all $s\in S$, we have that $P(S)=1$. For disjoint sets $B_1,\ldots, B_m$, using A2 we get additivity:
\[ P \left( \bigcup_{i=1}^{m}B_i \right) =\E\left(\sum_{i=1}^m I_{B_i}\right)= \sum_{i=1}^m P(B_i) .\]

Let us now show that $P(B)\geq 0$ for all $B\in \cF$.  Indeed, suppose that  $P(B)=-\eps$ for some $\eps>0$.
This implies (by A1 and A2) that  $\E (I_B +\eps)=0$; on the other hand, $I_B(s) + \eps >0$ for all $s\in S$, which contradicts both A3 and A3$'$,
and this ends the proof.

\eitem

\medskip

From now on, the quintuplet  $(S,\cF,\cS,\E,P)$ will be assumed to be as defined above in
the framework of axioms A1-A2-A3 and conditions $C1-C5$ (if not specified differently).
 The letter $X$ will be reserved for
elements of $\cS$, and $P$ is a set function derived from $\E$ as in  (\ref{aprob0}).

\medskip

\qitem{Lemma (positivity)}{\label{posit} \em If for all $s\in S$, $X(s)\geq 0$,  then $\E X \geq 0$.
}

\medskip

{\bf Proof. } Suppose   $X(s)\geq 0$  for all $s\in S$, and $\E X = -\eps$, for some $\eps >0$. Then by A1 and A2, $\E (X+\eps ) =0$, whereas $X+\eps >0$.
This contradicts A3. Therefore, if $X\geq 0$, then $\E X \geq 0$.

\eitem

\medskip

\qitem{Lemma}{\label{equ3} Assuming that A1 and A2 hold, the axioms A3 and A3$'$ are equivalent.}

{\bf Proof.} In  Lemma \ref{aprob} we already proved the property that $P$ is a FAP
follows with either  A3 or  A3$'$, so we may use that property in both parts of the  present proof.

Assume that A1-A2-A3 holds and suppose that for all $s\in S$, $X(s)\geq 0$  and $\E X =0$.
Then by A3 it follows that $P(X=0)=1$. Therefore, A3$'$ holds.

Now assume that A1-A2-A3$'$ hold, but not A3.  Then there exists $X\in \cS$ such that: a) $\E X =0$, $P(X=0)<1$, $X\geq 0$,
or b) $\E X =0$, $P(X=0)<1$, $X\leq 0$. In the case a), using A3$'$ we find that $P(X=0)=1$, which is a contradiction to $P(X=0)<1$.
The case b) can be reduced to a) with the function  $Y=-X$.

\medskip

Due to the equivalence established in  the Lemma \ref{equ3}, in the rest of article we refer to A3 as being either A3 or A3'.

\medskip

\qitem{Remark}{\label{hered} Suppose that we have a quintuplet $(S,\cF,\cS,\E,P)$, that satisfies assumptions C1-C5 and Axioms A1-A2,
 where $P$ is defined  with (\ref{aprob0}). Let $S^{\ast}\in \cF$ be a subset of $S$ with $P(S^{\ast})>0$, such that $X\cdot I_S^{\ast}\in \cS$.
 Define $\cF^{\ast}= \{ B\cap S^{\ast} \; |\; B\in \cF\}$.
 Next, for each $X\in \cS$ we define the  function
  $X^{\ast}:=X|_{S^{\ast}}$ - that is, $X$
restricted to $S^{\ast}$ and let $\cS^{\ast}$ be the space of all those  mappings.
 We define a linear functional $\E^{\ast}$ on $\cS$ as
  \[  \E^{\ast}(X^{\ast}) = \frac{1}{P(S^{\ast})}\E (X\cdot I_{S^{\ast}}), \qquad   X^{\ast}= X|_{S^{\ast}}.\]
  and the set function $P^{\ast}$ as
  \[ P^{\ast} (B^{\ast}) =\E^{\ast} (I_{B^{\ast}})= \frac{P(B^{\ast})}{P(S^{\ast})}=
  \frac{P(B\cap S^{\ast})}{P(S^{\ast})}, \qquad B^{\ast} = B\cap S^{\ast}\in \cF^{\ast}.\]

 In this way we get a new quintuplet with $(S^{\ast},\cF^{\ast},\cS^{\ast},\E^{\ast},P^{\ast})$,
and it is not difficult to see that the new quintuplet inherits  conditions $C1-C5$ and axioms A1-A2 , as well as $A3$ if it is satisfied
with the original  quintuplet.}

\eitem

In the next Lemma we prove Markov's inequality from the axioms.

\qitem{Lemma}{\label{lmarkovi}  Let $X\in \cS$ and $X(s)\geq 0$ for all $s\in S$. Then

\be
\label{markovi}
 P(X>\eps ) \leq \frac{\E X}{\eps}, \qquad \eps >0.
 \ee

}

{\bf Proof. } Since $X\geq 0$, we use C4 to conclude that

\[ \E X = \E (X\cdot I_{0\leq X\leq \eps}) + \E (X\cdot I_{X>\eps}).  \]
Further,we have
\[ X(s)\cdot I_{0\leq X\leq \eps}(s) \geq 0,\quad X(s)\cdot I_{X>\eps}(s) \geq \eps I_{X>\eps}(s), \]
and then use positivity (lemma \ref{posit}) and A1-A2 to conclude that $\E X \geq \eps P(X>\eps)$.

\eitem

\medskip

\qitem{Lemma}{\label{camis4} Assuming A1-A2, C1-C5 and Markov's inequality, A3 holds if  $P$ is countably additive probability.}

\medskip

{\bf Proof. }  Let $X \in \cS$ such that $X\geq 0$ and  $\E X=0$.  We need to show
that $P(X>0)=0$.
By Markov's inequality we have that $P(X>\eps) =0$ for any $\eps >0$, and so, using the countable additivity
we get
\[ P(X>0) = P\left( \bigcup_{n=1}^{\infty} \left\{ X > \frac{1}{n} \right\} \right) = \lim_{n\tends +\infty} P\left( X > \frac{1}{n}  \right)= 0, \]
as desired.

\eitem

\qitem{Remark}{\label{camok} Consider the  case where $P$ is a countably additive probability on $(S,\cF)$,
 $\cF$  is a sigma algebra, and  $\E X = \int_S X(s)\di P(s)$. Axioms A1-A2 and conditions C1-C5 are clearly satisfied, and Markov's inequality
 can be proved from properties of
 the integral, so by Lemma \ref{camis4}, the axiom A3 also holds.}\kraj

\eitem

Let us now recall some facts about FAPs. A probability $P$  which is defined on an  algebra $\cF$
of subsets of the set $S$, is {\em purely finitely additive} if $\nu \equiv 0$ is the only countably additive measure with the property that
$\nu (B) \leq P(B)$ for all $B\in \cF$. A purely finitely additive probability $P$ is {\em strongly finitely additive-SFAP} if there exist
countably many disjoint sets
$H_1,H_2,\ldots \in \cF$ such that
\be
\label{expart}
\bigcup_{i=1}^{+\infty} H_i=S \qquad  {\rm and}\qquad  P(H_i)=0\quad \mbox{for all $i$}.
\ee

For every probability $P$ on $\cF$ there exists a countably
additive probability $P_c$ and a purely finitely additive probability $P_d$ such that
$ P=\lambda P_c + (1-\lambda) P_d$ for some $\lambda \in [0,1]$.
This decomposition is unique (except  for $\lambda=0$ or $\lambda=1$, when it is trivially non-unique).

\medskip

\qitem{Lemma}{\em \label{sfpn4} Assuming axioms A1-A2, conditions C1-C5 and positivity, if $P$ is  a SFAP,
the condition of Axiom A3 is not satisfied.}

\medskip

{\bf Proof. } Let $P$ be a SFAP, and let $H_i$ be a partition of $S$ as in (\ref{expart}).
Define $X(s)=1/i$ if $s\in H_i$.  Further,
\[ X(s) \leq 1\cdot I_{H_1}(s) + \frac{1}{2} \cdot I_{H_2}(s) +\cdots + \frac{1}{k} \cdot I_{\{H_k \cup H_{k+1}\cup \cdots \}}\]
and so (by positivity) $0\leq \E X \leq 1/k$ for every $k>0$, hence $\E X =0$. This contradicts A3.

\eitem

\qitem{Example}{\label{ex1} Let $S= [0,+\infty)$ and let $P$ be the probability defined by the non-principal ultrafilter of Banach
limit as $s\tends +\infty$.
Let $X(s)=e^{-s}$. Then $X\geq 0$ and $\E X =0$, but $P(X=0)=0$.  In this case the convex hull $K(X) = (0,1]$ and $\E X \not\in K(X)$.
}

\eitem

\qitem{Theorem}{\label{ich} \em Let $(S,\cF,\cS,\E,P)$ be a quintuplet as defined above, and let
 $X = (X_1,\ldots, X_n)$, where $X_i \in \cS$ for all $i$. Assuming that axioms A1-A2-A3 and conditions C1--C5 hold,  $\E X$
  belongs to the
convex hull of the set
$X(S)=\{ X(t)\; |\; t\in S\}\subset \R^n$.}

{\bf Proof. } Without loss of generality we may assume that for all $i$, $\E X_i=0$ (otherwise if $\E X_i =c_i$ we
can observe $\E (X_i-c_i) =0$). Let $K$ denotes the convex hull of the set
$X(S)\in \R^n$. We now prove that $0\in K$ by induction on $n$. Let $n=1$. By A3, $\E X=0$ implies that either
$X(s)=0$ for some $s\in S$ or there are $s_1,s_2 \in S$ such that $X(s_1)>0$ and $X(s_2)<0$.
In both cases it follows that $0\in K$.

Now assume that the statement of the theorem is valid for  all dimensions from $1$ to $n-1$ for all quintuplets  $(S,\cF,\cS,\E,P)$
 that satisfy the conditions mentioned in the statement of the theorem. Let now $X$ be a vector function with values in $\R^n$.

If every  hyperplane $\pi$ that contains $0$  has the property that
the set $X(S)$  has an non-empty intersection with both of two open half-spaces with $\pi$ as boundary, then $0\in K$ (see \cite{rockaf} for details).
 Otherwise, suppose that
 \be
\label{lz}
 L(s):=\sum_{k=1}^n a_k X_k(s)\geq 0 \quad \mbox{for every $s\in S$},
 \ee
 with some  real numbers $a_1,\ldots, a_n$, such that $\sum a^2_k >0$. By linearity (A2) we have that $\E L(t)=0$, which is (A3$'$)
 possible together with (\ref{lz}) only if $L(t)=0$ for all $t\in S\setminus N$, where $\mu (N)=0$. Assuming that $a_n \neq 0$, we find that
  \be
 \label{xnt}
 X_n(s)= -\sum_{k=1}^{n-1} \frac{a_k}{a_n} X_k(s) \quad\mbox{for every $s\in S\setminus N$}.
 \ee

 In other words,  a separating hyperplane exists  only if there exists a linear relation among $n$ given functions with probability one. In order to
 eliminate $X_n$ and to reduce the system to $n-1$ functions, we
  consider  functions $X^{\ast}_i(s)=X_i(s)$ on the restricted domain $S^{\ast}=S\setminus N$, ($i=1,\ldots, n-1$) and the corresponding
  functional $\E^{\ast} $.
 Let $K^{\ast}$ be convex hull of $X^{\ast} (S^{\ast})\in \R^{n-1}$

By hereditary property (Remark \ref{hered}),  we have that $\E^{\ast} (X_i^{\ast}) = \E (X_i\cdot I_{S\setminus N})=0$, by C4.
Note that
$K^{\ast}\subset K$. By induction assumption, the statement of
the theorem holds for dimension $n-1$, and so
\be
 \label{summ}
 \sum_{i=1}^m \lambda_i X_k(t_i) =0 \qquad k=1,\ldots, n-1,
 \ee
 with some $t_1,\ldots, t_m \in S^{\ast}\subset S$ (here we use the fact that $X^{\ast}_i(s) = X_i(s)$ for $s\in S^{\ast}$).
 Finally, using (\ref{xnt}) and  (\ref{summ}) we find that also
 \be
 \label{xnnn}
 \sum_{i=1}^m \lambda_i X_n(t_i) =0,
 \ee
and so, the statement of the theorem holds for dimension $n$.

\eitem

\qitem{Remark}{\label{suff} Theorem \ref{ich} provides sufficient conditions for $\E X$ to belong to the convex hull of $X(S)$. However, by inspection
of the proof, we can see that the Axiom 3 is also necessary, assuming $A1-A2$ and conditions $C1-C5$.

}

Now, as a corollary to Theorem \ref{ich} using Caratheodory's theorem on representation of convex hull in finite dimension,
we get the following result:

\qitem{Theorem}{\label{carad} \em Assume that axioms A1-A2-A3 and conditions C1--C5 hold on $(S,\cF,\cS,\E,P)$.
Let  $X = (X_1,\ldots, X_n)$, where $X_i \in \cS$ for all $i$. Then
there are points $t_0,\ldots, t_n$ and a discrete probability law
given by probabilities $\lambda_0,\ldots, \lambda_n$ so that
\be
\label{disc}
\E X_i = \sum_{j=0}^n \lambda_j X_i (t_j),\qquad i=1,\ldots n.
\ee

}
\eitem

The Theorem \ref{carad} is the most general mean value theorem for axiomatic integral. Due to Remark \ref{camok}, the statement of this theorem
applies with $\E X_i = \int_S X_i(s)\di \mu (s)$, $\mu$ is a countably additive probability measure on $(S,\cF)$ where
$\cF$ is a sigma algebra, and $X_i$ are $(S,\cF)-(\R, \cB)$ measurable and integrable functions ($\cB$ is Borel sigma field on $\R$).

In the next section we consider the case of continuous
functions $X_i$.

\sect{Mean value theorem for continuous multivariate mappings }

\qitem{Definition}{A path from a point $a$ to point $b$ in a topological space $S$ is a continuous mapping $f: \ [0,1]\rightarrow S$
such that $f(0)=a$ and $f(1)=b$. A space $S$  is path connected if for any two points $a,b \in S$ there exists a path that connects them.}

\eitem

Let us remark that any topological vector space is path connected. A path that connects points $a$ and $b$ is given by
$f(\lambda)=\lambda a + (1-\lambda)b$, $\lambda \in [0,1]$. The same is true for a convex subset $S$ of any topological space.

\medskip

The following result is a generalization of Theorem B.

\qitem{\bf Theorem}{\label{chullt}\em Let $X: t\mapsto X(t)$, $t\in S$, be a continuous function defined on a path connected topological space
$S$ with values  in $\R^n$, and let $K$ be the convex hull of the set $X(S)$.
Then each $v\in K$ can be represented as a convex combination of $n$ or fewer points of the set $X(S)$.}

\medskip

{\bf Proof. } By Carath\'eodory's theorem, any $v\in K$ can be represented as a
convex combination of at most $n+1$ points of the set $X(S)$.
 Without loss of generality, assume that $v=0$.
Therefore, there exist $t_j\in S$ and $v_j\geq 0$, $0\leq j\leq n$, such that $t_i\ne t_j$ for $i\ne j$,
 $v_0+\cdots + v_n=1$, and
 \be
 \label{vct1}
v_0 x(t_0)+v_1x(t_1)+\cdots + v_n x(t_n)=0.
 \ee
We may also assume that all $n+1$ points $x(t_j)$ do not belong to one hyperplane in $\R^n$ (in particular, $x(t_i)\neq x(t_j)$ for $i\neq j$)
and that the numbers $v_j$ are all
positive; otherwise, at least one term from (\ref{vct1}) can be  eliminated. Now we apply the following reasoning:

 Denote by $p_j(x)$, $0\leq j \leq n$,
the coordinates of the vector $x\in \R^n$ with respect to the coordinate system with the origin at $0$, and with the vector
base consisted of vectors $x(t_j)$, $j=1,\ldots, n$ (that is, $x=\sum_{j=1}^n p_j (x) x(t_j)$). Then from (\ref{vct1}) we find that
$p_j\left(x(t_0)\right) = -v_j/v_0<0$, $j=1,\ldots, n$, i.e. the coordinates of the vector
$x(t_0)$ are negative.
The coordinates of vectors $x(t_j)$, $j=1,2,\ldots, n$ are non-negative: $p_j\left( x(t_j)\right) =1$ and
$p_k\left( x(t_j)\right)=0$ for $k\ne j$. Now consider a path $t=t(\lambda)$, $\lambda \in [0,1]$ which connects
points $t_0$ and $t_1$, so that $ t(0) = t_0$ and $t(1)=t_1$. The functions $\lambda \mapsto p_j(x(t(\lambda))):=f_j(\lambda)$
are continuous as mappings from $[0,1]$ to $\R^n$ and $f_j(0) <0$ for all $j=1,\ldots,n$, whereas $f_1(1)=1$ and $f_j(1)=0$ for $j>1$.
Therefore, for each of functions $f_j$ there exists one or more points $\lambda\in (0,1]$ such that $f_j(\lambda)=0$.
Since the set $N=\cup_{j=1}^n f_j^{-1}(\{0\})$ is closed and non-empty subset of $(0,1]$, there exists $\lambda_0 =\min S $, $\lambda_0 >0$. Let
$\bar{t}:=t(\lambda_0)$. From this construction it follows that there exists (at least one) $k$ such that $p_k\left( x(\bar{t})\right)=0$
and $p_j\left( x(\bar{t})\right) <0$ for $j\neq k$. Hence,

\[x(\bar{t}) - \sum_{1\leq j\leq n, j\neq k} p_j\left( x(\bar{t})\right) x(t_j)=0 \]

 and it follows that  $v=0$ is a convex combination of  points $x(\bar{t})$ and $x(t_j)$, $j=1,\ldots, n$, $j\ne k$.

\eitem

As a direct corollary to Theorems \ref{ich} and  \ref{chullt}, we have the following mean value theorem for continuous multivariate mappings.

\qitem{Theorem}{\label{mainth}\em Let $(S,\cF,\cS,\E,P)$ be a quintuplet as defined in Section 2, where $S$ is a
 path connected topological space. Under conditions C1-C5 and axioms A1-A3,  let $X_i, i=1,\ldots, n$ be continuous functions from $\cS$.
  Then there exist points $t_1,\ldots, t_n$ in $S$,
and non-negative numbers $\lambda_1,
\ldots, \lambda_n$, with $\sum_{i=1}^n \lambda_i = 1$, such that

\[ \E X_k = \sum_{i=1}^n \lambda_i X_k(t_i),\qquad k=1,\ldots, n.\]

}

\eitem

\qitem{Remark}{\label{remval} Since Theorem \ref{chullt} is independent of axioms and conditions of Section 2, the statement of Theorem \ref{mainth}
holds whenever the Theorem \ref{carad} holds. In particular, it holds whenever $\E X_i = \int_S X_i(s)\di \mu (s)$, where $\mu$ is a countably additive
probability measure.

In fact, all what Theorem \ref{mainth} says is that we can save one point in the representation
(\ref{disc}) of Theorem \ref{carad}. Although it might look not much significant, in some applications it makes difference. For example, Karamata's
representation for covariance in \cite{karamceb} based on Kowalewski's original result for $n=2$, strongly depend on two points
and nothing similar can be derived with three points.
}

\newpage

\sect{Some particular cases and open problems}

\medskip

\qitem{Riemann and Lebesgue integral on $\R^d$}{\label{RL} For $d=1$, let $-\infty<a<b<+\infty$ and let $X$ be a Riemann integrable function
 on $[a,b]$. Define
 \be
 \label{rint}
 \E X = \frac{1}{b-a} \int_a^b X(s)\di s .
 \ee
In terms of previous notations, $S=[a,b]$, and $\cS$ is the family of Riemann integrable functions on $\cS$. A natural choice for
algebra $\cF$ of sets related to Riemann integral should be the algebra of intervals in $[a,b]$, which can be defined as the collection of all
subintervals of $[a,b]$ (including singletons and empty set) and their finite unions.
The corresponding probability is defined then as follows: If $A=\cup_{i=1}^k J_i$ where $J_i$ are intervals,
\[ P(A) =  \frac{1}{b-a} \int_a^b I_{A}(s)\di s = \sum_{i=1}^k \lambda(J_i), \]
where $\lambda(J_i)$ is the length of $J_i$.  This is Jordan probability measure on $[a,b]$, and it is well known that, even for continuous bounded
functions, the set $\{s\in S: X(s)\leq c\}$ where $c$ is a real number,  does not obligatory belong to $\cF$ (this is probably
first shown by an example in \cite{frin33}). This fact makes it impossible to use our system of axioms directly, because  condition C3 does not hold.
Nevertheless, we can proceed by noticing that the algebra $\cF$ as described above is a sub-algebra of the Borel sigma algebra $\cB$, generated by open
sets in $[a,b]$, and since Riemann and Lebesgue integrals coincide if the integrand is Riemann integrable, we can proceed it this way.
}
In more common notations, let us consider a general case of a functional $\E $ based on Lebesgue integral on $\R^d$, $d\geq 1$:
\be
\label{leint}
\E f = \frac{1}{V(D)} \int_D f(s_1,\ldots, s_d)\di s_1\ldots \di s_d,
\ee
or, in shorthand,
\[ \E f =\frac{1}{V(D)} \int_D f(s)\di \lambda(s), \]
where $\lambda$ is the Lebesgue measure on $\R^d$ restricted to $D$, and $V(D) = \lambda(D)>0$, where $D$ is a convex (or in more generality, path connected)
 subset of $\R^d$.  The underlying probability measure in our
construction of Section 2 is $P(\cdot)= \frac{1}{V(D)} \lambda (\cdot)$. Let $f_1,\ldots, f_n$ be continuous functions $D\mapsto \R$
such that $\E f_i$ as defined  in (\ref{leint}). Then, by Theorem \ref{mainth}, we have that there are points $x_1,\ldots, x_n \in D$ and
non-negative numbers $\lambda_1,\ldots, \lambda_n$ with $\sum_{i=1}^n \lambda_i =1$ such that

\beqn
\frac{1}{V(D)} \int_D f_1(s_1,\ldots, s_d)\di s_1\ldots \di s_d &=& \lambda_1 f_1(x_1) +\cdots \lambda_n f_1(x_n) \\
            \vdots                                              &=& \vdots  \\
\frac{1}{V(D)} \int_D f_n(s_1,\ldots, s_d)\di s_1\ldots \di s_d &=& \lambda_1 f_n(x_1) +\cdots \lambda_n f_n(x_n) \\
\eeqn

In words, this result shows that for  an arbitrary system of $n$ integrals with continuous integrands, there exists an
exact quadrature rule with $n$ points in $D$, with coefficients $\lambda_i$ which are the same for all integrals.
Note that $x_i$ are $d$-dimensional points, so in fact here we have $dn$ scalar parameters.

\eitem

\qitem{Integrals with respect to countably additive probability measure}{As already noted, Theorem \ref{mainth} holds for all integrals based
on countably additive probability measure. Suppose that $S$ is a path connected topological space, and let $\mu$ be a  countably additive probability  measure
on $(S,\cF)$, where $\cF$ is a sigma algebra of Borel subsets of $S$. Let $f_1,\ldots, f_n$ be continuous mappings from $S$ to $\R$, and
suppose that
\be
\label{exrp}
 \E f_i= \int_S f(s)\di\mu(s),\qquad i=1,\ldots,n
 \ee
is finite. In particular, let $S=C[0,T]$, the space of continuous functions on the interval $[0,T]$ with the supremum norm. Then the measure $\mu$ defines a
stochastic process on $[0,T]$, $s(t), 0\leq t\leq T$  and $f_i(s)$ is a continuous functional of trajectories of the process. By Theorem \ref{mainth}, we have
that the system of expectations (\ref{exrp}) can be represented as
\[ \E f_i= \sum_{j=1}^n \lambda_j f_i (x_j) ,\qquad i=1,\ldots,n \]
for some $x_j \in C[0,T]$ and $\lambda_j\geq 0$ with $\sum_j \lambda_j=1$.

}

\eitem

\qitem{Open question}{We showed in Section 2 that Axiom 3 does not hold for (integrals based on) SFAPs,
so by Remark \ref{suff}, the
mean value theorem does not  hold in this case. It would be of interest to describe  classes of finitely additive probabilities
for which Axiom 3 holds or does not hold, in terms of some structural properties of measures.}

\eitem

\medskip

{\bf Acknowledgment.} This work is supported by grants 11402 and III 44006 from Ministry of Education, Science and Technological Development of
Republic of Serbia.

\bigskip

\noindent{University of Belgrade}

\noindent{Faculty of Electrical Engineering}

\noindent{Department of Applied Mathematics}

\noindent{Bulevar Kralja Aleksandra 73}

\noindent{11120 Belgrade, Serbia}

\noindent{\verb+emerkle@etf.rs+}

\end{document}